\documentclass[11pt]{article}
\usepackage{graphicx}
\usepackage{amssymb,amsfonts,amsthm}
\usepackage{enumerate}
\usepackage[T2B]{fontenc}
\usepackage[cp1251]{inputenc}

\usepackage{enumitem}
\usepackage{amsmath,amsthm,amssymb}
\usepackage{amscd}
\usepackage{amsfonts}
\usepackage{color}

\newcommand{\N}{{\mathbb N}}
\newcommand{\R}{{\mathbb R}}

\newtheorem{theorem}{Theorem}

\theoremstyle{definition}

\usepackage{imakeidx}
\makeindex[name=g,title=Subject index]
\usepackage[columns=1]{idxlayout}
\begin{document}
\renewcommand{\theequation}{\thesection.\arabic{equation}}
\bigskip

\title{The isoperimetric spectrum of finitely presented groups}

 \author{M.V. Sapir\thanks{Supported in part by the NSF grant DMS-1500180.}}

\date{}
\maketitle

\begin{abstract}
The isoperimeric spectrum consists of all real positive numbers $\alpha$ such that $O(n^\alpha)$ is the Dehn function of a finitely presented group. In this note we show how a recent result of Olshanskii completes the description of the isoperimetric spectrum modulo the celebrated Computer Science conjecture (and one of the seven Millennium Problems) $\mathbf{P=NP}$ and even a formally weaker conjecture.
\end{abstract}

The goal of this note is to show that the recent paper by  Olshanskii \cite{Ol} completes a description of the isoperimetric spectrum of finitely presented groups modulo the $\mathbf{P=NP}$ conjecture.

Since in this note we consider only polynomially bounded functions $\N\to \R$, we call two functions $f,g$ {\em equivalent} if $af(n)\le g(n)\le bf(b)$ for some positive constants $a,b$. 
 
Recall that Brady and Bridson \cite{Bridson} called the set of all real numbers $\alpha\ge 1$ such that $n^{\alpha}$ is equivalent to the Dehn function of a finitely presented group the {\em isoperimetric spectrum}. When it was introduced, it was known only that all natural numbers belong to the isoperimetric spectrum (the free nilpotent group of class $c$ with at least 2 generators has Dehn function $n^{c+1}$ \cite{B}), and that by Gromov's theorem the intersection of the isoperimetric spectrum with the open interval $(1,2)$  is empty. It is obvious also that the isoperimetric spectrum is a countable set since the set of all finite group presentation is countable. Bridson \cite{Bridson} found the first examples of non-integral numbers in the spectrum.

Soon after, we proved in \cite{SBR} that for  $\alpha\ge 4$ to be in the isoperimetric spectrum, it is enough that $\alpha$ is computed in time $\le 2^{2^m}$.
 Recall \cite{SBR} that a real number $\alpha$ is called computable in time $T(m)$ where $T(m)$ is a function $\N\to \N$, if there exists a
deterministic Turing machine which, given a natural number $n$,
computes a binary rational approximation of $\alpha$ with an error
at most $1/2^{n+1}$ in at most $T(n)$ steps.
Thus all algebraic numbers $\ge 4$ and many transcendental numbers such as $\pi+1$ are in the isoperimetric spectrum. On the other hand, we proved in \cite{SBR} that every number in the isoperimetric spectrum can be computed in time $\le 2^{2^{2^{cm}}}$ for some constant $c$. It can be seen from the proof of this result that the number of $2$'s in this estimate can be reduced to two as in the lower bound if we had $\mathbf{P=NP}$. I mentioned it (without a proof) in my ICM talk \cite{SapirICM} and my Bulletin of Mathematical  Sciences survey \cite{SapirBMS}.

For $\alpha\le 4$, the situation was more complicated. On the one hand the tools used in \cite{SBR} were too weak to handle $\alpha\le 4$. On the other hand, Brady, Bridson, Forester  and Shankar found more numbers from the interval $(2,4)$ in the isoperimetric spectrum, showing that the set of these numbers is dense in the interval $(2,4)$ \cite{BradyBridson} and even contains all rational numbers \cite{B2}. Their numbers from the isioperimetric spectrum were constructed using algebraic rather than computational properties. (Note also that the groups constructed in \cite{Bridson, BradyBridson, Brady, B2}, are given by very small presentations comparing to  the groups in \cite{SBR} and are subgroups of CAT(0) groups which is quite remarkable.) But the paper by Olshanskii \cite{Ol} showed that the intersection of the isoperimetric spectrum with $(2,4)$ can be described in the same terms as in \cite{SBR}.

Combining results of \cite{SBR} and \cite{Ol} we get:

\begin{theorem} [The first part of Corollary 1.4 from \cite{SBR}, and Corollary 1.4  from \cite{Ol}]\label{1} If \\ a number $\alpha\ge 2$ can be computed in time $\le 2^{c2^m}$ for some $c$, then $\alpha$ belongs to the isoperimetric spectrum.
\end{theorem}

Notice that one needs to modify a little the proof of the first part of Corollary 1.4 from \cite{SBR} to obtain the estimate $\le 2^{c2^m}$ instead of $\le 2^{2^m}$. For this, one should take natural number $d>\log_2 c$ which is a power $2^k$ for some $k$ and consider $d$-ary representations of numbers instead of binary representations as in \cite{SBR}. Each $d$-ary digit of $\alpha$ is $k$ binary digits. So if the first $m$ binary digits of $\alpha$ are computed in time $\le 2^{c2^m}$, then the first $m$ $d$-ary digits of $\alpha$ are computed in time at most $d^{d^m}$ and the rest of the proof of \cite[Corollary 1.4]{SBR} carries by replacing $2$ by $d$ everywhere.

Now we will prove the main result of the note. 

\begin{theorem} \label{2} Provided $\mathbf{P=NP}$, a number $\alpha$ is in the isoperimetric spectrum if and only if it can be computed in time $\le 2^{c2^m}$ for some $c\ge1$.
\end{theorem}

\proof Theorem \ref{1} gives one part of Theorem \ref{2}. 

To prove the other part, suppose that $n^\alpha$ is the Dehn function of a finitely presented group. Then
by \cite[Theorem 1.1]{SBR}, $n^\alpha$ is equivalent to the time function $T(n)$ of
some (non-deterministic) Turing machine $M$ which recognizes the word problem in a finitely presented group. 

As explained by Emil Je\v r\'abek \cite{MO}, the following property follows from $\mathbf{P=NP}$:

\begin{quote}
(*) There is a deterministic Turing machine $M'$ computing a function $T'(n)$ which is equivalent to $T(n)$ and having time function at most $T(n)^d$ for some constant $d$.
\end{quote}

Since $\alpha\ge2$, for any $n>0$
$$\epsilon_1 n^\alpha\le T'(n) \le \epsilon_2n^\alpha$$
for some positive constants $\epsilon_1\le 1$ and $\epsilon_2\ge 1$.
Let number $n_0$ be such that $2^n>\log_2 (\epsilon_2/\epsilon_1)$ for every $n\ge n_0$.
Let $q=[\alpha]+1$.

Given this machine $M'$ with, say, $k$ tapes, consider the following Turing machine $M''$ which will
calculate the first $m$ digits of $\alpha$ (for every $m$). This machine has $k+3$ tapes
with tape $k+3$ being the input tape. It starts with number $m$ in binary
written on tape
$k+3$ and all other tapes empty. Then it calculates the number $n=2^{2^{m+n_0}}$
and
writes it on tape $k+1$ (using tape $k+2$ as an auxiliary tape and
cleaninig it after $n$ is computed). Then $M''$ turns on the machine $M'$
and produces $T'(n)$ on tape $k+2$. Then it calculates $p=[(\log_2T'(n)-\log_2\epsilon_1)/2^{n_0}]$
and writes it on tape $k+2$. Notice that $\alpha\log_2 n+\log_2\epsilon_1\le
\log_2 T'(n)\le  \alpha\log_2 n+\log_2 \epsilon_2$. Therefore
$$[\alpha 2^m]\le p\le
\alpha 2^m+\log_2(\epsilon_2/\epsilon_1)/2^{n_0}.$$
Hence $p=[\alpha 2^m]$, so $p/2^m$ a rational approximation of $\alpha$ which is
within $1/2^m$ from $\alpha$.
From the construction of $M''$, it is clear that the time complexity
of $M''$ does not exceed $2^{c2^m}$
for some constant $c$.
\endproof

As explained by Emil Je\v r\'abek \cite{MO}, Property (*) follows also from the property $\mathbf{E=\Sigma^E_2}$, which is not known to imply $\mathbf{P = NP}$.  Here $\mathbf{E}$ denotes the class of languages recognized by deterministic Turing machines in time $2^{O(n)}$ and $\mathbf{\Sigma^E_2}$ is the second level of the exponential hierarchy (with linear exponent). Thus in Theorem \ref{2}, one can replace $\mathbf{P=NP}$ by $\mathbf{E=\Sigma^E_2}$.

{\bf Acknowledgement.} I am grateful to Emil Je\v r\'abek for answering my question \cite{MO} and several useful remarks. I am also grateful to Jean-Camille Birget and A. Yu. Olshanskii for their comments.


\begin{thebibliography}{99}

\bibitem{B} G. Baumslag, C. F. Miller III, and H. Short, Isoperimetric inequalities and the homology of groups, Invent. Math. 113 (1993), 531--560.


\bibitem{BradyBridson} N. Brady and M. Bridson, There is only one gap in the isoperimetric spectrum, Geometric and Functional Analysis, 10 (2000), 1053--1070.

\bibitem{B2} Noel Brady, Martin R. Bridson, Max Forester, Krishnan Shankar,
Snowflake groups, Perron-Frobenius eigenvalues and isoperimetric spectra. Geom. Topol. 13 (2009), no. 1, 141--187. 


\bibitem{Brady} Noel Brady, Max Forester,
Snowflake geometry in CAT(0) groups.
J. Topol. 10 (2017), no. 4, 883--920. 



\bibitem{Bridson}  M. Bridson, Fractional isoperimetric inequalities and subgroup distortion. 
J. Amer. Math. Soc. 12 (1999), no. 4, 1103--1118.

\bibitem{MO} Mathoverflow, question 307629. 

\bibitem{Ol} A. Yu. Olshanskii, Polynomially-bounded Dehn functions of groups,  arXiv:1710.00550, accepted in Journal of Combinatorial Algebra, 2018. 

\bibitem {SapirICM} M. V.  Sapir,
Algorithmic and asymptotic properties of groups.  International Congress of Mathematicians. Vol. II, 223--244, Eur. Math. Soc., Z\"urich, 2006.

\bibitem{SapirBMS} M. V. Sapir, Asymptotic invariants, complexity of groups and related problems. Bull. Math. Sci. 1 (2011), no. 2, 277--364.

\bibitem {SBR}  M. V. Sapir, J. C. Birget, E. Rips,
Isoperimetric and isodiametric functions of groups,
Annals of Mathematics, 157, 2(2002), 345-466.

\end{thebibliography}
\end{document}